\documentclass[paper=a4,english,fontsize=11pt,parskip=half,abstract=true]{scrartcl}
\usepackage{babel}
\usepackage[utf8]{inputenc}
\usepackage[T1]{fontenc}
\usepackage[left=20mm,right=20mm,top=30mm,bottom=30mm]{geometry}
\usepackage{amsmath}
\usepackage{amsthm}
\usepackage{amssymb}
\usepackage[shortlabels]{enumitem}
\usepackage{thmtools}
\usepackage{mathtools}
\usepackage{embedfile}
\mathtoolsset{centercolon} 
\usepackage[bookmarks=true,
            pdftitle={Counterexamples to congruences in finite p-groups},
            pdfauthor={Benjamin Sambale},
            pdfkeywords={p-groups, products of conjugacy classes, products of characters},
            pdfstartview={FitH}]{hyperref}

\newtheorem{Thm}{Theorem} 
\numberwithin{equation}{section}

\allowdisplaybreaks[1]

\renewcommand{\phi}{\varphi}

\newcommand{\FF}{\mathbb{F}}

\newcommand{\Irr}{\operatorname{Irr}}

\title{Counterexamples to congruences in finite $p$-groups}
\author{Benjamin Sambale\footnote{Institut für Algebra, Zahlentheorie und Diskrete Mathematik, Leibniz Universität Hannover, Welfengarten 1, 30167 Hannover, Germany,
\href{mailto:sambale@math.uni-hannover.de}{sambale@math.uni-hannover.de}}}
\date{\today}

\begin{document}
\frenchspacing
\maketitle
\embedfile[desc={GAP code verifying the counterexamples},
           mimetype={text/plain}]{kkinv.g}
\begin{abstract}\noindent
Let $K$ be a conjugacy class of a finite $p$-group $G$. Héthelyi, Külshammer and the author asked whether the number of conjugacy classes contained in $KK^{-1}$ is always congruent to $1$ modulo $p-1$, and whether, dually, the number of irreducible constituents of $\chi\overline{\chi}$ is congruent to $1$ modulo $p-1$ for every $\chi\in\Irr(G)$. We construct a group of order $7^7$ for which both congruences fail.
\end{abstract}


Let $G$ be a finite $p$-group. For a conjugacy class $K$ of $G$ let $K^{-1}:=\{a^{-1}:a\in K\}$, and let $\eta(K)$ be the number of conjugacy classes of $G$ contained in $KK^{-1}$. Since all conjugacy classes have $p$-power length, $\eta(K)\equiv|KK^{-1}|\pmod{p-1}$. For a complex character $\psi$ of $G$ let $\Irr(\psi)$ be the set of its irreducible constituents. Motivated by sharp lower bounds for $\eta(K)$ and $\lvert\Irr(\chi\overline{\chi})\rvert$ due to Adan-Bante~\cite{AdanBanteCl,AdanBanteChar}, we asked in \cite{HKSconj} whether
\[\eta(K)\equiv 1\pmod{p-1}\tag{P1}\]
for every conjugacy class $K$ of $G$, and whether
\[\lvert\Irr(\chi\overline{\chi})\rvert\equiv 1\pmod{p-1}\tag{Q1}\]
for every $\chi\in\Irr(G)$. There, (P1) and (Q1) were proved for $p\leq 3$, for $|G|\leq p^6$, for groups of nilpotency class at most $2$ and for metacyclic groups. We show that both congruences fail in general.

\begin{Thm}\label{main}
Let $R:=\FF_7[s]/(s^5)$, let $P:=\langle 1+s,\,1+s^2+s^4\rangle\leq R^\times$, and let $G:=R\rtimes P$, where $P$ acts on the additive group of $R$ by multiplication. Then $|G|=7^7$. 
\begin{enumerate}[(a)]
    \item\label{maineta} For the conjugacy class $K$ of $1\in R$ we have $\eta(K)=33\equiv 3\pmod 6$. 
    \item\label{mainirr} Let $\lambda\in\Irr(R)$ be given by $\lambda(a):=\zeta^{\tau(a)}$, where $\zeta:=e^{2\pi i/7}$ and $\tau(a)\in\FF_7$ is the coefficient of $s^4$ in $a$. Then $\chi:=\lambda^G\in\Irr(G)$ and $\lvert\Irr(\chi\overline{\chi})\rvert=81\equiv 3\pmod 6$.
\end{enumerate}
In particular, (P1) and (Q1) fail for $G$.
\end{Thm}

\begin{proof}
Since $(1+s)^7=1+s^7=1$ and $(1+s^2+s^4)^7=1$, both generators of $P$ have order $7$. The coefficient of $s$ in $(1+s)^i$ is $i$, so $\langle 1+s\rangle\cap(1+s^2R)=1$. Hence $|P|=49$ and $|G|=7^5\cdot 7^2=7^7$.

For $0\ne a\in R$ let $\nu(a)$ be the largest $k$ with $a\in s^kR$. For $a,b\ne 0$ we have $ab\ne 0$ whenever $\nu(a)+\nu(b)\leq 4$. Every $r\in P\setminus\{1\}$ has the form $r=(1+s)^i(1+s^2+s^4)^j$, and $\nu(r-1)=1$ if $7\nmid i$ and $\nu(r-1)=2$ otherwise.

Let $P-P:=\{a-b:a,b\in P\}$. We claim that $P$ acts freely on $(P-P)\setminus\{0\}$ by multiplication. Let $a,b,c\in P$ with $a\ne b$ and $c(a-b)=a-b$. Since $a-b=b(ab^{-1}-1)$ and $b$ is a unit, we get $(c-1)(ab^{-1}-1)=0$. As $\nu(ab^{-1}-1)\leq 2$, the previous paragraph forces $c=1$. A computer calculation shows $|P-P|=1569$. Thus $(P-P)\setminus\{0\}$ splits into $1568/49=32$ orbits of length $49$.

\ref{maineta} Since $R$ is abelian, the conjugacy class in $G$ of an element $a\in R$ is its orbit $Pa$. In particular $K=P$, $KK^{-1}=P-P$, and the conjugacy classes in $KK^{-1}$ are the $P$-orbits on $P-P$. Hence $\eta(K)=1+32=33$.

\ref{mainirr} The map $\tau\colon R\to\FF_7$ is a homomorphism of additive groups. For $r\in R$ let $\lambda_r(a):=\zeta^{\tau(ra)}$. Then $\lambda_r\in\Irr(R)$, and $(R,+)\to\Irr(R)$, $r\mapsto\lambda_r$ is a homomorphism. It is injective, because for $r\ne 0$ we can write $r=s^ku$ with a unit $u$ and get $\tau(r\cdot s^{4-k}u^{-1})=1$. Since $|R|=\lvert\Irr(R)\rvert$, it is an isomorphism, and $\lambda=\lambda_1$. For $g\in P$ we have $\lambda_r^g(a)=\lambda_r(ga)=\lambda_{gr}(a)$, so the action of $P$ on $\Irr(R)$ corresponds to multiplication on $R$. In particular, the stabilizer of $\lambda$ in $P$ is trivial, and $\chi=\lambda^G$ is irreducible by \cite[Theorem~6.11]{Isaacs}. Moreover $\chi_R=\sum_{g\in P}\lambda_g$ and $\overline{\chi}=\overline{\lambda}^G$, so \cite[Problem~5.3]{Isaacs} yields
\[\chi\overline{\chi}=(\chi_R\overline{\lambda})^G=\sum_{g\in P}(\lambda_{g-1})^G.\]
Here $(\lambda_0)^G=(1_R)^G$ is the sum of the $49$ linear characters of $G/R\cong P$. For $g\ne 1$ the element $g-1$ lies in $(P-P)\setminus\{0\}$, so $\lambda_{g-1}$ has trivial stabilizer and $(\lambda_{g-1})^G$ is irreducible of degree $49$. By Clifford's theorem \cite[Theorem~6.2]{Isaacs}, $(\lambda_{d})^G=(\lambda_{d'})^G$ holds if and only if $d$ and $d'$ lie in the same $P$-orbit. Since every element of $P-P$ lies in the $P$-orbit of some $g-1$ with $g\in P$, all $32$ orbits on $(P-P)\setminus\{0\}$ occur. Therefore $\lvert\Irr(\chi\overline{\chi})\rvert=49+32=81$.
\end{proof}

Analogous constructions yield counterexamples of order $p^7$ for all primes $p\ge 7$.
There are no counterexamples of order $5^7$, but there is a counterexample of order $5^{12}$ with $R:=\FF_5[s,t]/(s^4,t^2)$ and $P:=\langle 1+s,\ 1+t,\ 1+s^2,\ 1+st+s^3\rangle\leq R^\times$. 

The group $G$ of \autoref{main} is \texttt{SmallGroup}$(7^7,78341)$ in the Small Groups Library. The file attached to this pdf verifies the counterexamples of order $7^7$ and $5^{12}$ in GAP~\cite{GAP4161}.

\section*{Acknowledgment}
The work \cite{HKSconj} is based on my diploma thesis written in 2008. At the time, I found a counterexample to a stronger conjecture, called (P3) in \cite{HKSconj}, of order $3^7$ (note that (P1) and (Q1) hold for $p=3$ in general). The groups of order $7^7$ were not available at the time. Claude Opus~5.5 (Anthropic) took 30 minutes to find the present counterexample without brute-forcing the Small Groups Library. I have checked all the details and take full responsibility.

\end{document}